\documentclass[11pt]{amsbook}

\title[A priori estimates for the Yamabe problem]{A priori estimates for the Yamabe problem in the non-locally conformally flat case}

\author{Fernando C. Marques}

\newtheorem{thm}{Theorem}[section]

\newtheorem{prop}[thm]{Proposition}

\newtheorem{lem}[thm]{Lemma}

\newtheorem{definition}[thm]{Definition}

\newtheorem{cor}[thm]{Corollary}

\numberwithin{equation}{section}

\def\g{\rightarrow}
\def\R{\mathbb{R}}
\def\Rn{{\mathbb{R}}^n}
\def\d{\partial}

\begin{document}

\renewcommand{\qedsymbol}{{\bf q.e.d.}}

\begin{abstract}
Given a compact Riemannian manifold $(M^n,g)$, with positive Yamabe quotient, not conformally
diffeomorphic to the standard sphere, we prove 
a priori estimates for solutions to the Yamabe problem. We restrict ourselves to the dimensions where the Positive  
Mass Theorem is known to be true, that is, when $n \leq 7$. We also show that, when $n \geq 6$, the Weyl tensor has to 
vanish at a point where  solutions to the Yamabe equation blow up.
\end{abstract}

\maketitle

\tableofcontents


\section{Introduction}

Let $(M^n,g)$ be a compact Riemannian n-manifold, $n \geq 3$, without boundary. The classical Yamabe problem
consists of finding a metric $\tilde{g}$, conformally related to $g$, with constant scalar curvature on $M$.
It can be considered as a generalization of the classical Uniformization Theorem on Riemann surfaces
 to the setting of higher dimensional manifolds.

In analytical terms, the problem is equivalent to show the existence of a positive solution $u$ to the equation
\begin{equation}\label{eq}
\Delta_gu -\frac{n-2}{4(n-1)}R_gu + Ku^\frac{n+2}{n-2} = 0 \ {\rm on }\ M,
\end{equation}
 where $\Delta_g$ denotes the Laplace-Beltrami operator associated with the metric $g$, $R_g$
denotes the scalar curvature of the metric $g$, and $K$ is a constant. The linear operator
$L_g = \Delta_g -\frac{n-2}{4(n-1)}R_g$ is called the conformal Laplacian of the metric $g$.

The  solution of the Yamabe problem was an outstanding achievement since, for the first time, it was
given a very satisfactory existence theory to a  
nonlinear partial differential
equation involving a critical exponent.
 After the initial paper in 1960 by Yamabe \cite{YAMABE60}, which contained an error, contributions made by 
Trudinger \cite{TRUDINGER68},
Aubin \cite{AUBIN76}, and finally by Schoen \cite{SCHOEN84} in 1984,  solved the problem completely in the affirmative. 

In this paper we shall be interested in the set of solutions to the Yamabe problem. 
When the first eigenvalue of the 
conformal Laplacian $L_g$ is negative, it is not difficult to see
that the solution is unique. If this eigenvalue is zero, 
the equation becomes linear and then solutions are unique up to a constant. Therefore the only interesting case 
left is the positive one.

We know that when the underlying manifold is the sphere $(\mathbb{S}^n, g_0)$, endowed with the standard metric, M. Obata's theorem (see \cite{OBATA71}) classifies
all solutions to the equation, and this set 
is noncompact in the $C^2$ topology. On the other hand, the standard sphere is the only compact manifold with a noncompact group of conformal diffeomorphisms, so
one should expect different behavior in the other cases.

In \cite{SCHOEN91}, R. Schoen proved the compactness, in the $C^2$ topology, of the set of solutions to the Yamabe equation, in the 
positive case, for every locally conformally flat manifold not conformally diffeomorphic to the sphere. He also suggested a strategy to prove these  a priori $C^{2,\alpha}$ estimates in the non-locally conformally flat case, based on Pohozaev-type identities. In \cite{LIZHU99}, Y. Y. Li and M. Zhu followed these lines and proved the theorem in dimension
3, in which case standard estimates on the blowing up solutions are sufficient for a Pohozaev identity to be applied. 
The compactness result in dimensions 4 and 5 was obtained by O. Druet in \cite{DRUET2003,DRUET2004}.
Other compactness theorems for the Yamabe equation, in the locally conformally flat case, are proved in 
\cite{HANLI99} for manifolds with boundary, and in \cite{POLLACK93} for singular solutions on the sphere. 

Our main result in this paper is the following a priori estimates theorem in the general case: 
\begin{thm}\label{comp}
 Let $(M^n,g)$ be a smooth closed Riemannian manifold with positive Yamabe quotient, 
not conformally equivalent to $(\mathbb{S}^n, g_0)$. Assume $n \leq 7$. Then, 
for every $\epsilon > 0$,  there exists a 
positive constant $C = C(\epsilon, g)$ so that
$$
\left\{
\begin{array}{lr}
   &1/C \leq u \leq C \ {\rm and} \\
   &||u||_{C^{ 2,\alpha}(M)} \leq C 
\end{array}
\right.
$$
for every $u \in \cup_{1+\epsilon \leq p \leq \frac{n+2}{n-2}}\mathbb{M}_p$, where $0 < \alpha < 1$ and
$$
\mathbb{M}_p = \{u > 0:\Delta_gu -\frac{n-2}{4(n-1)}R_gu + Ku^p = 0 \ {\rm on }\ M\}.
$$
\end{thm}

These estimates clearly imply the compactness of the set of solutions to the Yamabe equation in the $C^2$ topology.  We restrict ourselves to the 
dimensions
covered by the Positive Mass Theorem due to Schoen and Yau \cite{SCHOENYAU79}, i.e., $n \leq 7$, since the final global argument in our proof depends essentially on this result. 
Due to technical problems concerning singularities of minimizing hypersurfaces, the Positive Mass Theorem is still not known to be true for dimensions greater than 7. 

Our  result will follow from a contradiction between a local restriction coming from a Pohozaev-type identity and a 
global argument provided by the Positive Mass Theorem.  In order to accomplish that we will need a careful blowup analysis of solutions, part of it inspired
by the work of C. C. Chen and C. S. Lin \cite{CHENLIN98}. We introduce some 
new symmetry estimates which allow us to control how close the blowing up solutions get to some specific rotationally
symmetric functions. These symmetric functions will be  solutions to the corresponding critical ($p=\frac{n+2}{n-2}$) equation in $\Rn$.  Since we are also dealing with subcritical equations, estimates on $\tau = \frac{n+2}{n-2}-p$ are also given. 

We would like
to point out that the symmetry estimates (Proposition \ref{simetria}), when $n \leq 5$, are just as in the conformally flat case. However, when $n \geq 6$ we can no longer expect these same estimates to hold. This is   because, in general, the asymptotic expansion of the Green function for the conformal Laplacian has additional terms (see \cite{LEEPARKER87}). 


One important difficulty we must overcome when $n \geq 6$, pointed out by Schoen in \cite{SCHOEN91}, is to show conformal flatness of the
metric to a sufficiently high order at a blowup point. This is needed in order to apply  the Positive Mass Theorem when $n=6$ or $7$. That is the content of our next theorem, where $W_g$ denotes the Weyl tensor of the metric $g$:

\begin{thm}\label{weylintrod}
Assume $n \geq 6$ and let $u_i$ be  a sequence of positive solutions to (\ref{eq}).
Suppose $x_i \g \overline{x}$ is a sequence of points such that $u_i(x_i) \g \infty$ as $i \g \infty$.
Then
$$
W_g(\overline{x}) = 0.
$$
\end{thm}

In general one should expect, as indicated by Schoen in   \cite{SCHOEN90}, that at a blowup point $\overline{x}$ we must have :
$$
\nabla^kW_g(\overline{x}) = 0
$$ for all $0 \leq k \leq \frac{n-6}{2}$. The proof of Theorem \ref{weylintrod} also relies upon
the symmetry estimates and a Pohozaev-type identity. 


It is not difficult to check that the  Theorem \ref{comp} implies the existence of a solution to the Yamabe problem. This is because solutions $u_i$ to the subcritical equations, with $p_i \g \frac{n+2}{n-2}$ as $i \g \infty$, can be constructed by standard variational methods. Actually one can say more. Another  consequence of the compactness theorem is the computation of the total Leray-Schauder
degree of all solutions to equation (\ref{eq}): 

\begin{thm}
Let $(M^n,g)$ be a smooth closed Riemannian manifold with positive Yamabe quotient, 
not conformally equivalent to $(\mathbb{S}^n, g_0)$, $4 \leq n \leq 7$.
Then, if $\Lambda$ is sufficiently large,
$$
{\rm deg}(F,\Omega_{\Lambda},0)=-1,
$$
where $F(u) = u + L_g^{-1}(E(u)u^\frac{n+2}{n-2})$, $E(u) = -\int_M uL_g(u) dv_g$ and 
$$
\Omega_{\Lambda} = \{u \in C^{2,\alpha}(M): \min_{M} u > \Lambda^{-1}, ||u||_{2,\alpha} < \Lambda\}.
$$
\end{thm}
This theorem follows from Theorem \ref{comp} and arguments given by Schoen in \cite{SCHOEN91}.

The author was just recently communicated that Y. Y. Li and L. Zhang  have independently proved the 
same compactness result.


{\bf Acknowledgements.}
The content of this paper is  part of the author's doctoral thesis \cite{MARQUES2003}. The author  would like to dedicate this paper to the memory of  his friend and advisor Prof. Jos\'{e} F. Escobar. His encouragement was invaluable for the completion of this work. While the author was at Cornell University, he was fully supported by CNPq-Brazil.


\section{A Pohozaev-type identity}

In this section we will establish a Pohozaev-type identity which will be very useful in the subsequent  blowup
analysis. 

Suppose 
$u:B_{\rho}(0) \setminus \{0\} \subset \Rn \g \R $ 
is a  positive $C^2$ solution to the equation 
\begin{equation}\label{equacao}
a^{ij}(x)\d_{ij}u + b^i(x)\d_iu + c(x)u + K(x)u^p = 0,
\end{equation}
where $p \neq -1$, $K \in C^1$ and  $a^{ij},b^i,c$ are continuous functions, $1 \leq i,j \leq n$. Here we are using the 
summation convention. 

Define 
\begin{eqnarray}\label{invpoh}
&& \ \ P(r,u) = \\
&&\int_{|x|=r}(\frac{n-2}{2}u\frac{\d u}{\d r}
-\frac{|x|}{2}|\nabla u|^2 + |x||\frac{\d u}{\d r}|^2 + 
\frac{1}{p+1}K(x)|x|u^{p+1})d\sigma(r), \nonumber
\end{eqnarray}
whenever $0 < r < \rho$.

The following lemma gives the Pohozaev-type identity we are interested in.

\begin{lem}
Given $0 < s \leq r < \rho$,
\begin{eqnarray*}
&&P(r,u)-P(s,u)= 
-\int_{s \leq |x| \leq r}(x^k\d_ku + \frac{n-2}{2}u)
A(u)dx  \\
&& \ \ \ \ \ \ \ \ \ \ \ +\frac{1}{p+1}\int_{s \leq |x| \leq r}(x^k\d_kK(x))u^{p+1} dx \\
&& +(\frac{n}{p+1}-\frac{n-2}{2})
\int_{s \leq |x| \leq r}K(x)u^{p+1}dx,
\end{eqnarray*}
where $A(u) = (a^{ij}-\delta_{ij})\d_{ij}u + b^i\d_iu + cu$.
\end{lem}

\begin{proof}
Multiplying the equation (\ref{equacao}) by $x^k\d_ku$, and integrating
over the set $\{x:s \leq |x| \leq r\}$, we obtain 
\begin{equation}\label{poh1}
\int_{s \leq |x| \leq r}(x^k\d_ku)(\Delta u + A(u) + K(x)u^p)dx = 0.
\end{equation}

Integration by parts gives:
\begin{eqnarray*}\label{poh2}
&&\int_{s \leq |x| \leq r}(x^k\d_ku)\d_{ii}u dx = \nonumber\\
&&- \int_{s \leq |x| \leq r}(\delta_i^k\d_ku\d_iu 
+ \frac12 x^k\d_k[(\d_iu)^2])dx \nonumber\\
&&\ \ \ \ \ \ +\frac{1}{r}\int_{|x|=r}(x^k\d_ku)(x_i\d_iu) - \frac{1}{s}\int_{|x|=s}(x^k\d_ku)(x_i\d_iu) = \nonumber\\
&&- \int_{s \leq |x| \leq r}(\delta_i^k\d_ku\d_iu -\frac{n}{2}(\d_iu)^2)dx 
-\frac{r}{2}\int_{|x|=r}(\d_iu)^2
+\frac{s}{2}\int_{|x|=s}(\d_iu)^2 \nonumber\\
&&\ \ \ \ \ \ +\frac{1}{r}\int_{|x|=r}(x^k\d_ku)(x_i\d_iu) - \frac{1}{s}\int_{|x|=s}(x^k\d_ku)(x_i\d_iu),
\end{eqnarray*}
and summing  over $i=1,\dots,n$, we obtain
\begin{eqnarray}\label{poh3}
&&\int_{s \leq |x| \leq r}(x^k\d_ku)\Delta u dx =\nonumber\\
&&\frac{n-2}{2}\int_{s \leq |x| \leq r}|\nabla u|^2 dx 
-\frac{r}{2}\int_{|x|=r}|\nabla u|^2+\frac{s}{2}\int_{|x|=s}|\nabla u|^2\nonumber\\
&&\ \ \ \ \ \ \ + r\int_{|x|=r}(\frac{\d u}{\d r})^2 - s\int_{|x|=s}(\frac{\d u}{\d r})^2.
\end{eqnarray}

Also
\begin{eqnarray}\label{poh4}
&&\int_{s \leq |x| \leq r}(x^k\d_ku)K(x)u^p dx = \nonumber\\
&&\ \ \frac{1}{p+1}\int_{s \leq |x| \leq r}x^k\d_k(u^{p+1})K(x)dx=\nonumber\\
&&\ \  \ -\frac{n}{p+1}\int_{s \leq |x| \leq r}K(x)u^{p+1}dx - 
\frac{1}{p+1}\int_{s \leq |x| \leq r}(x^k\d_kK(x))u^{p+1} dx \nonumber\\
&&\ \ \ \ \ \ +\frac{r}{p+1}\int_{|x|=r}K(x)u^{p+1} - \frac{s}{p+1}\int_{|x|=s}K(x)u^{p+1}.
\end{eqnarray}

On the other hand, multiplying the equation (\ref{equacao}) by u,
and once again integrating by parts, we get 
\begin{eqnarray}\label{poh5}
&&\int_{s \leq |x| \leq r}|\nabla u|^2dx =\nonumber \\
&&\ \ \ \int_{s \leq |x| \leq r}(uA(u)+K(x)u^{p+1})dx + 
\int_{|x|=r}u\frac{\d u}{\d r} - \int_{|x|=s}u\frac{\d u}{\d r}.
\end{eqnarray}

 Now we substitute equalities (\ref{poh3}), (\ref{poh4}) 
and equality (\ref{poh5}) in equality (\ref{poh1}) and the Pohozaev identity 
follows by rearranging terms.
\end{proof}

When u is a solution to the equation (\ref{equacao}) in the entire ball, by taking the limit as $s \g 0$ we get
\begin{eqnarray}\label{pohozaev}
&&\ \ P(r,u)= 
-\int_{|x| \leq r}(x^k\d_ku + \frac{n-2}{2}u)A(u)dx  \\
 &&\ \ \ \ \ \ +\frac{1}{p+1}\int_{|x| \leq r}(x^k\d_kK(x))u^{p+1} dx +\nonumber \\
&&(\frac{n}{p+1}-\frac{n-2}{2})\int_{|x| \leq r}K(x)u^{p+1}dx. \nonumber
\end{eqnarray}

Integrating by parts once more, we can also get
\begin{eqnarray}\label{pohozaev2}
&&P(r,u)=-\int_{|x| \leq r}(x^k\d_ku + \frac{n-2}{2}u)
((a^{ij}-\delta^{ij})\d_{ij}u + b^i\d_iu)dx \\
&&\ \ \ \ \ \ +\int_{|x| \leq r}(\frac12 x^k\d_kc + c)u^2 dx
-\frac{r}{2}\int_{|x|=r}cu^2 d\sigma(r) \nonumber \\
&&\ \ \ +\frac{1}{p+1}\int_{|x| \leq r}(x^k\d_kK(x))u^{p+1} dx \nonumber \\
&&+(\frac{n}{p+1}-\frac{n-2}{2})\int_{|x| \leq r}K(x)u^{p+1}dx.\nonumber
\end{eqnarray}


\section{Conformal scalar curvature equation}\label{confinvariance}

In this section we will introduce the partial differential equation we are interested in, and  we shall discuss
some of its properties related to conformal deformation of metrics.    

Let $\Omega \in \Rn$ be an open set, and suppose $g$ is a 
Riemannian metric in $\Omega$. Suppose also $f$ is a positive $C^1$ function defined in $\Omega$. 

Consider a positive  $C^2$ function $u$  satisfying
\begin{equation}\label{eqbasica}
\Delta_{g}u -c(n)R_{g}u + Kf^{-\tau}u^{p} = 0 \ {\rm in } \  \Omega,
\end{equation}
where $c(n) = \frac{n-2}{4(n-1)}$, $K = n(n-2)$,  
$1 <p\leq \frac{n+2}{n-2}$ and $\tau = \frac{n+2}{n-2} - p$. We will use the notation $R_{g}$ for the scalar curvature
of $g$. The operator $L_{g} = \Delta_{g}-c(n)R_{g}$ is called the 
{\it conformal Laplacian} of the metric $g$. 

When $p = \frac{n+2}{n-2}$, this partial differential equation is intimately related to conformal geometry, particularly 
when one studies conformally related metrics with constant scalar curvature. More specifically, given a positive solution $u$,
the metric $u^\frac{4}{n-2}g$ has constant scalar curvature equal to $4n(n-1)$.

Now let us describe an important feature of solutions to that type of  equation. Let $u$ be a solution to equation (\ref{eqbasica}) and
choose $\overline{x} \in \Omega$. Given $s > 0$, define the renormalized function
$$
v(y) = s^\frac{2}{p-1}u({\rm exp \ }_{\overline{x}}(sy)).
$$

Then
$$
L_{h}v + K\tilde{f}^{-\tau}v^{p} = 0,
$$
where $\tilde{f}(y)=f(sy)$ and the components of the metric $h$ in normal coordinates are given by $h_{kl}(y)=g_{kl}(sy)$.
The important point here is that $v$ satisfies an equation of the same type.

The equation is also conformally invariant in the following sense. Suppose $\tilde{g} = \phi^\frac{4}{n-2}g$ is a metric conformal to $g$. Let us recall
\begin{equation}\label{conflapla}
L_{\tilde{g}}(\phi^{-1}u) = \phi^{-\frac{n+2}{n-2}}L_g(u)
\end{equation}
for any function $u$, and 
\begin{equation}\label{confcurv}
R_{\tilde{g}}=-c(n)^{-1}\phi^{-\frac{n+2}{n-2}}L_g(\phi).
\end{equation}
(See \cite{LEEPARKER87}).

 Therefore, if $u$ is a solution to (\ref{eqbasica}), then $\phi^{-1}u$ satisfies
$$
L_{\tilde{g}}(\phi^{-1}u) + K(\phi f)^{-\tau}(\phi^{-1}u)^p = 0,
$$
which is again an equation of the same type. 

This will have very important consequences in what follows. We will study  sequences of solutions 
$u_i$ to equation  (\ref{eqbasica}). When what we want to study is conformally invariant, we are allowed to
 replace  $u_i$ by another sequence of functions $v_i = \phi_i^{-1}u_i$, at the same time replacing the metric $g_i$ by
$\tilde{g}_i = \phi_i^\frac{4}{n-2}g_i$, as long as we have a uniform control on the conformal factors $\phi_i$. In this paper, there will be two examples of such a procedure.

First, we can suppose the metric $g_i$ has positive scalar curvature in a small ball centered at some
fixed point $x_i$. To see this, fix $\sigma > 0$ small and let $\phi$ be the first eigenfunction of $\Delta_{g}$ with respect to 
the Dirichlet condition:
\begin{equation}
\left\{
\begin{array}{lr}
&\Delta_g \phi + \lambda_1 \phi = 0 \ {\rm in }\  B_{2\sigma}(x) \\
&\phi = 0 \ {\rm on }\  \d B_{2\sigma}(x).
\end{array}
\right.
\end{equation}

Recall that the corresponding eigenspace is one-dimensional and we can choose $\phi > 0$ on $B_{2\sigma}(x)$. Now, since $\lambda_1 \g \infty$ as $\sigma \g 0$, we can choose $\sigma$ small
enough so that
$$
\Delta_g \phi - c(n)R_{g}\phi < 0
$$
in $B_{\sigma}(x)$. Defining $\tilde{g} = \phi^\frac{4}{n-2} g$, relation (\ref{confcurv}) implies
 $R_{\tilde{g}} > 0$ in $B_{\sigma}(x)$. Moreover, if we take as a conformal factor the solution $\psi$ to:
\begin{equation}
\left\{
\begin{array}{lr}
&\Delta_g \psi - c(n)R_g \psi = 0 \ {\rm in }\  B_{\sigma}(x) \\
&\psi = 1 \ {\rm on }\  \d B_{\sigma}(x),
\end{array}
\right.
\end{equation}
we can also have $R_{\tilde{g}} = 0$. 

The second example is related to the so-called conformal normal coordinates. (See \cite{LEEPARKER87}).  Given an integer $N \geq 2$, 
there exists a positive function $\phi$ (which can be constructed explicitly), such that, setting $\tilde{g}=\phi^\frac{4}{n-2}g$,   the volume element satisfies:
$$
\det(\tilde{g}_{ij}) = 1 + O(r^N),
$$
in $\tilde{g}$-normal coordinates around $x$, where $r = d_{\tilde{g}}(x,\cdot)$. 
This allows us to simplify the local asymptotic analysis. For example, in conformal normal coordinates around $x$, 
 $R_{\tilde{g}} = O(r^2)$ and $\Delta R_{\tilde{g}}(x) = -\frac16 |W_{\tilde{g}}(x)|^2$, where $W$ stands for the Weyl tensor.


\section{Isolated and isolated simple blowup points}

In this section we will define isolated and isolated simple blowup points and we shall
discuss their basic properties. The results in this section are well-known 
in the locally conformally flat setting (\cite{LI95}, \cite{SCHOENZHANG96}) and in general when $n=3$ (\cite{LIZHU99}). We will
slightly modify their proofs in \cite{LIZHU99} to make them work in any dimension.

Let $\Omega \subset \Rn$ be an open set, and suppose $g_i$ is a sequence of
Riemannian metrics in $\Omega$ converging, in the $C^2_{loc}$ topology, to a metric
g. Suppose also that $f_i$ is a sequence of positive $C^1$ functions converging in the $C^1_{loc}$ topology
to a positive function $f$. 

We will consider a sequence $u_i$ of  positive $C^2$ functions satisfying
\begin{equation}\label{isol1}
L_{g_i}u_i + Kf_i^{-\tau_i}u_i^{p_i} = 0 \ {\rm in }\  \Omega,
\end{equation}
where $c(n) = \frac{n-2}{4(n-1)}$, $K = n(n-2)$,  
$1+\epsilon_0 <p_i\leq \frac{n+2}{n-2}$ for some $\epsilon_0 > 0$ and $\tau_i = \frac{n+2}{n-2} - p_i$. 

We will sometimes omit the subscript $i$, for the sake of simplicity, and we will use the symbols $c,C$ to denote various
positive constants.


\begin{definition}\label{defisol}
We say that $\overline{x} \in \Omega$ is an {\rm isolated blowup point} for $u_i$ if 
there exists a sequence $x_i \in \Omega$, converging to $\overline{x}$, so that:
\begin{enumerate}
\item $x_i$ is a local maximum point of $u_i$;
\item $M_i := u_i(x_i) \g \infty$ as $i \g \infty$;
\item there exist $r,C>0$ such that
\end{enumerate}
\begin{equation}\label{ineqisol}
u_i(x) \leq Cd_{g_i}(x,x_i)^{-\frac{2}{p_i-1}}
\end{equation}
for every $x \in B_r(x_i)\subset \Omega$. Here $B_r(x_i)$ denotes the geodesic ball of radius $r$, centered at $x_i$,
with respect to the metric $g_i$.
\end{definition}

{\bf Remark:} In various parts of the text, we will identify $x_i$ with the origin, that meaning we are making use of normal
coordinates in a small ball around $x_i$. More precisely, we will sometimes write $u_i(x)$ instead of
$u_i({\rm exp}_{x_i}(x))$ and $|x|$ instead of $d_{g_i}(x,x_i)$, and those functions will be defined in balls centered at 0.

Note that  the definition of isolated blowup points  is invariant under renormalization,
which was described  in the last section. This follows from the fact that, if $v(y) = s^\frac{2}{p_i-1}u(sy)$, then
$$
u(x) \leq C|x|^{-\frac{2}{p_i-1}} \Leftrightarrow v(y) \leq C|y|^{-\frac{2}{p_i-1}}.
$$


The first result concerning isolated blowup points is the following
Harnack inequality.

\begin{prop}\label{harnack}
Suppose that $u_i$ is a sequence of positive functions satisfying equation (\ref{isol1}) and
assume $x_i \g \overline{x}$ is an isolated blowup point. Then there exists
a constant $C > 0$ such that
$$
\max_{\frac{s}{2}\leq d_{g_i}(x,x_i) \leq 2s} u_i(x) \leq
C\min_{\frac{s}{2}\leq d_{g_i}(x,x_i) \leq 2s} u_i(x),
$$
where $0 < s < \frac{r}{3}$.
\end{prop}

\begin{proof}

Let $(x^1,\dots,x^n)$ be normal coordinates with respect to the metric $g_i$ on the ball $B_r(x_i)$. (See  remark
after Definition \ref{defisol}).

Define 
$$
v_i(y) = s^{\frac{2}{p_i-1}}u_i(sy), 
$$
where $|y| < 3$.

Then, as discussed in the last section,
$$
L_{h_i}v_i(y) + K\tilde{f_i}^{-\tau_i}v_i^{p_i}(y) = 0,
$$
where $\tilde{f}_i(y)=f_i(sy)$ and $(h_i)_{kl}(y)=(g_i)_{kl}(sy)$, and we also know that
$$
v_i(y) \leq C|y|^{-\frac{2}{p_i-1}},
$$
whenever $|y| < 3$.

It follows from this last inequality that $v_i$ is uniformly bounded in
compact subsets of $B_3(0) \setminus \{0\}$. The Harnack inequality
for elliptic linear equations then implies that there exists $C>0$ such that
$$
\max_{\frac{1}{2}\leq |y| \leq 2} v_i(y) \leq
C\min_{\frac{1}{2}\leq |y| \leq 2} v_i(y).
$$
The result now follows directly.
\end{proof}

The Proposition \ref{harnack} clearly implies the so-called {\it spherical Harnack inequality} for isolated
blowup points. Namely, 
given $0 < s \leq \frac23 r$, there exists a positive constant $C$, not depending on $s$, such that 
\begin{equation}\label{sphericalHarnack}
\max_{d_{g_i}(x,x_i)=s} u_i(x) \leq
C\min_{d_{g_i}(x,x_i)=s} u_i(x).
\end{equation}


Define $U_0(y) = (1+|y|^2)^{\frac{2-n}{2}}$. It is not difficult to check that
$$
\Delta U_0(y) + KU_0^\frac{n+2}{n-2}(y) = 0.
$$

The next proposition says that, in the case of an isolated blowup point, the functions $u_i$, when renormalized, converge in
the $C^2$ topology to the rotationally symmetric function $U_0$.

\begin{prop}\label{propisol}
Let $u_i$ be a sequence of positive functions satisfying the equation (\ref{isol1}) and
$x_i \g \overline{x}$ be an isolated blowup point. Assume that $R_i \g \infty$ and
$\epsilon_i \g 0$ are given. Then $p_i \g \frac{n+2}{n-2}$ and, after possibly passing to a subsequence,
\begin{equation}
||M_i^{-1}u_i(M_i^{-\frac{p_i-1}{2}}y)
-U_0(y)||_{C^2(B_{R_i}(0))} \leq \epsilon_i
\end{equation}
and
$$
\frac{R_i}{\log M_i} \g 0 \ {\rm as }\  i \g \infty.
$$
Here $M_i$ is as in Definition \ref{defisol}.
\end{prop}

\begin{proof}
Let $(x^1,\dots,x^n)$ be normal coordinates with respect to the metric $g_i$ on the ball $B_r(x_i)$. (See first remark
after Definition \ref{defisol}).

Define
$$
v_i(y) = M_i^{-1}u_i(M_i^{-\frac{p_i-1}{2}}y)
$$
for $|y| < rM_i^{\frac{p_i-1}{2}}$. Here $r$ is as in Definition \ref{defisol}.

Then
$$
L_{h_i}v_i(y) + K\tilde{f_i}^{-\tau_i}v_i^{p_i}(y) = 0,
$$
where $\tilde{f}_i(y)=f_i(M_i^{-\frac{p_i-1}{2}}y)$ and
 $(h_i)_{kl}(y)=(g_i)_{kl}(M_i^{-\frac{p_i-1}{2}}y)$.

Note also that
\begin{equation}\label{isol2}
  \left\{
     \begin{array}{lr}
    &v_i(0) = 1, \nabla v_i(0) = 0 \\
    &0 < v_i(y) \leq C|y|^{-\frac{2}{p_i-1}} \ {\rm for }\  |y| < rM_i^{\frac{p_i-1}{2}}.
  \end{array}
   \right.
\end{equation}

{\bf Claim:}
There exists $C>0$ such that $v_i(y) \leq C$, whenever 
$|y| <  rM_i^{\frac{p_i-1}{2}}$.

{\it Proof of  Claim:}
From properties (\ref{isol2}), we get
\begin{equation}\label{isol3}
v_i(y) \leq C,
\end{equation}
if $1 \leq |y| \leq rM_i^{\frac{p_i-1}{2}}$.

Now, from what was discussed in the previous section, up to a conformal deformation we can suppose our metrics have zero scalar curvature in small balls.  In particular their conformal Laplacians will satisfy the maximum principle. This implies
there exists $C > 0$ so that
$$
\min_{|y| \leq r}v_i(y) \geq C^{-1}\min_{|y|=r}v_i(y) \ \  \forall \ i,
$$
and $0<r \leq 1$. The spherical Harnack inequality (\ref{sphericalHarnack}) implies
\begin{eqnarray}\label{isol4}
\max_{|y| = r}v_i(y) &\leq& C\min_{|y| = r}v_i(y) \leq
C\min_{|y| \leq r}v_i(y)\nonumber\\
 &\leq& Cv_i(0) = C,
\end{eqnarray}
for $0<r \leq 1$. This and inequality (\ref{isol3}) imply the claim.

Standard elliptic estimates now imply that, after passing to a subsequence,
$v_i \g v > 0$ in $C^2_{loc}(\Rn)$, where
\begin{equation}
   \left\{
    \begin{array}{lr}
    &\Delta v(y) + Kv^p(y)=0 ,y \in \Rn \\
    &v(0)=1, \ \  
    \nabla v(0)=0,
   \end{array}
   \right.
\end{equation}
where $p = \lim_{i \g \infty}p_i$.
Here $\Delta$ denotes the Euclidean Laplacian.

A well-known theorem by Caffarelli, Gidas and Spruck \cite{CAFFAGIDASSPRUCK89}, states that we necessarily have $p = \frac{n+2}{n-2}$ and $v(y) = U_0(y)$.

The Proposition now follows easily.
\end{proof}


Now let us introduce the notion of an isolated simple blowup point.

Suppose $u_i$ is a sequence of positive functions satisfying equation (\ref{isol1}) and
$x_i \g \overline{x}$ is an isolated blowup point. Define
$$
\overline{u}_i(r)=
\frac{1}{\sigma_{n-1}r^{n-1}}\int_{\d B_r(x_i)}u_i d\sigma(r),
$$
where $\sigma_{n-1}$ denotes the area of a unit sphere in $\Rn$.
We are using  $g_i$-normal coordinates and integrating with respect to the Euclidean volume form.
\begin{definition}
 We say $x_i \g \overline{x}$ is an
{\rm isolated simple blowup point} if there exists a real number $0 < \rho < r$
such that the functions
$$
\hat{u}_i(r) = r^\frac{2}{p_i-1}\overline{u}_i(r)
$$  
have exactly one critical point in the interval $(0,\rho)$, for $i$ large.
\end{definition}

It is not difficult to see that Proposition \ref{propisol} implies that $\hat{u}_i$ has exactly 
one critical point in the interval $(0,R_iu_i(x_i)^{-\frac{p_i-1}{2}})$.
Moreover its derivative is negative right after the critical point.
As a result, if the blowup is isolated simple, then
$$
\hat{u}_i'(r) < 0
$$
for all $R_iM_i^{-\frac{p_i-1}{2}} \leq r < \rho$.


Now we turn to the first estimate on isolated simple blowup points.

\begin{prop}\label{propsimples}
Let $u_i$ be a sequence of positive functions satisfying equation (\ref{isol1}) and
$x_i \g \overline{x}$ be an isolated simple blowup point. Then there exists
a constant $C > 0$ and $0 < \rho_1 < \rho$ such that, for each $i$,
\begin{equation}\label{estsimples}
M_iu_i(x) \leq Cd_{g_i}(x,x_i)^{2-n},
\end{equation}
whenever $d_{g_i}(x,x_i) \leq \rho_1$. Moreover, if 
$R_iM_i^{-\frac{p_i-1}{2}} \leq d_{g_i}(x,x_i) \leq \rho_1$, then
\begin{equation}\label{estbaixo}
M_iu_i(x) \geq C^{-1} G_i(x_i,x),
\end{equation}
where $G_{i}$ is the Green function of $L_{g_i}$ with respect to the Dirichlet boundary condition on 
$B_{\rho_1}(x_i)$.
\end{prop}

\begin{proof}
We first need a slightly different estimate.

Let $\delta > 0$, and define $\lambda_i = (n-2-\delta)\frac{p_i-1}{2}-1$.

Let us  apply the Proposition \ref{propisol} to some $R_i \g \infty$
and $0 < \epsilon_i < e^{-R_i}$.

{\bf Claim 1:}
If $\delta$ is sufficiently small, there exist constants $0 <\rho_1 < \rho$ and
 $C>0$ such that 
\begin{eqnarray}
M_i^{\lambda_i} u_i(x) &\leq& C d(x,x_i)^{2-n+\delta}, \\
M_i^{\lambda_i}|\nabla u_i(x)| &\leq&  Cd(x,x_i)^{1-n+\delta},   \\
M_i^{\lambda_i} |\nabla^2 u_i(x)|     &\leq&  Cd(x,x_i)^{-n+\delta}, 
\end{eqnarray}
for every $x$ so that $R_iM_i^{-\frac{p_i-1}{2}} \leq d(x,x_i) \leq \rho_1$.

The proof of Claim 1 is analogous to the proof of the Lemma 3.3 in \cite{LIZHU99}.

{\bf Remark:}
It is not difficult to see that the previous estimates imply
\begin{equation}\label{estvi}
\left\{
\begin{array}{lr}
v_i(y) &\leq CM_i^{\delta \frac{p_i-1}{2}}(1+|y|)^{2-n},\\
|\nabla v_i(y)| &\leq CM_i^{\delta \frac{p_i-1}{2}}(1+|y|)^{1-n},\\
|\nabla^2 v_i(y)| &\leq CM_i^{\delta \frac{p_i-1}{2}}(1+|y|)^{-n}
\end{array}
\right.
\end{equation}
for any $|y| \leq \rho_1 M_i^\frac{p_i-1}{2}$.

Let us now estimate $\tau_i$.

{\bf Claim 2:}
There exists $C > 0$ such that
$$
\tau_i \leq 
\left\{
\begin{array}{lr}
&CM_i^{(-1+\delta)\frac{4}{n-2}+o(1)} \ {\rm if }\  n>4,\\
&CM_i^{(-1+\delta)2+o(1)}\log M_i \ {\rm if }\ n=4,
\end{array}
\right.
$$ 
 and, in particular, $M_i^{\tau_i} \g 1$ as 
$i \g \infty$. 

{\it Proof of Claim 2:}
We will apply the Pohozaev identity (\ref{pohozaev}) from Section 2 to $u_i$ on the ball of radius 
$\frac{\rho_1}{2}$:
\begin{eqnarray}\label{pohclaim2}
&&P(\frac{\rho_1}{2},u_i)= 
-\int_{|x| \leq \frac{\rho_1}{2}}(x^m\d_mu_i + \frac{n-2}{2}u_i)A_i(x)dx \nonumber \\
&&\ \ \ \ \ \ \ + (\frac{n}{p_i+1}-\frac{n-2}{2})\int_{|x|\leq \frac{\rho_1}{2}}
Kf_i^{-\tau_i}u_i^{p_i+1}dx \nonumber\\
&&-\frac{\tau_i}{p_i+1}\int_{|x|\leq \frac{\rho_1}{2}}Kf_i^{-\tau_i-1}(x^m\d_mf_i)u_i^{p_i+1}dx,
\end{eqnarray}
where
\begin{eqnarray*}
&&A_i(x) = (g^{kl}-\delta^{kl})(x)\d_{kl}u_i(x) \\
&&\ \ \ \ +(\d_kg^{kl}+ |g|^{-\frac12}\d_k(|g|^{\frac12})g^{kl})(x)\d_lu_i(x)
-c(n)R_g(x)u_i(x)
\end{eqnarray*}
and recall
\begin{eqnarray*}
&&P(\frac{\rho_1}{2},u_i)=
\int_{|x| = \frac{\rho_1}{2}}
(\frac{n-2}{2}u_i\frac{\d u_i}{\d r}-\frac{|x|}{2}|\nabla u_i|^2 
+ |x||\frac{\d u_i}{\d r}|^2) d\sigma \\
&&\ \ \ \ +\frac{1}{p_i+1}\int_{|x| = \frac{\rho_1}{2}}Kf_i^{-\tau_i}|x|u_i^{p_i+1}d\sigma.
\end{eqnarray*}

From Claim 1  we get
\begin{equation}\label{estP1}
|P(\frac{\rho_1}{2},u_i)| \leq cM_i^{-2\lambda_i}.
\end{equation}

Define
\begin{eqnarray}
&&\hat{A}_i(y) = 
(g^{kl}-\delta^{kl})(M_i^{-\frac{p_i-1}{2}}y)\d_{kl}v_i \nonumber\\
&&\ \ \ \ +M_i^{-\frac{p_i-1}{2}}(\d_kg^{kl}+|g|^{-\frac12}\d_k(|g|^{\frac12})g^{kl})(M_i^{-\frac{p_i-1}{2}}y)\d_lv_i \nonumber\\
&& -c(n)M_i^{-(p_i-1)}R_g(M_i^{-\frac{p_i-1}{2}}y)v_i.
\end{eqnarray}

The change of variables $y = M^\frac{p_i-1}{2}x$, the inequalities (\ref{estvi}) and the fact that the metric is euclidean
up to first order in normal coordinates yield
\begin{eqnarray}\label{estoutra}
&&|\int_{|x| \leq \frac{\rho_1}{2}}(x^m\d_mu_i + \frac{n-2}{2}u_i)A_i(x)dx| \nonumber\\
&&= M_i^{p\frac{2-n}{2}}M_i^\frac{n+2}{2}
|\int_{|y| \leq \frac{\rho_1}{2}M^\frac{p_i-1}{2}}(y^m\d_mv_i + \frac{n-2}{2}v_i)\hat{A}_i(y)dy|\nonumber\\
&&\leq CM_i^{p\frac{2-n}{2}}M_i^\frac{n+2}{2}M_i^{-(p_i-1)}M_i^{\delta(p_i-1)}
\int_{|y| \leq \frac{\rho_1}{2}M^\frac{p_i-1}{2}}(1+|y|)^{4-2n}dy\nonumber\\
&&\leq 
\left\{
\begin{array}{lr}
&CM_i^{(-1+\delta)\frac{4}{n-2}+o(1)} \ {\rm if }\  n>4,\\
&CM_i^{(-1+\delta)2+o(1)}\log M_i \ {\rm if }\ n=4.
\end{array}
\right.
\end{eqnarray}

Therefore, from the inequalities (\ref{estP1}) and (\ref{estoutra}) and the identity (\ref{pohclaim2}),
\begin{eqnarray}\label{estoutra1}
&&\frac{(n-2)\tau_i}{2(p_i+1)}\int_{|x|\leq \frac{\rho_1}{2}}Kf_i^{-\tau_i}u_i(x)^{p_i+1}dx \nonumber\\
&&\ \ \ \ -\frac{\tau_i}{p_i+1}\int_{|x|\leq \frac{\rho_1}{2}}Kf_i^{-\tau_i-1}(x^m\d_mf_i)u_i^{p_i+1}dx\nonumber\\
&&\leq 
\left\{
\begin{array}{lr}
&CM_i^{(-1+\delta)\frac{4}{n-2}+o(1)} \ {\rm if }\ n>4,\\
&CM_i^{(-1+\delta)2+o(1)}\log M_i \ {\rm if }\ n=4.
\end{array}
\right.
\end{eqnarray}
Since, from the Proposition \ref{propisol},
$$
\int_{|x|\leq R_iM_i^{-\frac{p_i-1}{2}}}u_i(x)^{p_i+1}dx \geq c > 0,
$$
we conclude that if we choose $\rho_1$ sufficiently small, then
\begin{eqnarray}
&&\frac{n-2}{2}\int_{|x|\leq \frac{\rho_1}{2}}Kf_i^{-\tau_i}u_i(x)^{p_i+1}dx \nonumber\\
&&- \int_{|x|\leq \frac{\rho_1}{2}}Kf_i^{-\tau_i-1}(x^m\d_mf_i)u_i^{p_i+1}dx \geq c > 0.
\end{eqnarray}
The result follows immediately from the inequalities (\ref{estoutra1}).

{\bf Claim 3:}
Given a small $\sigma > 0$, there exists a constant $C > 0$ such that
$$
\int_{B_\sigma(x_i)}u_i^{p_i} (x) dx \leq C M_i^{-1}.
$$

{\it Proof of Claim 3:}
Set $s_i=R_iM_i^{-\frac{p_i-1}{2}}$.

First note that, changing variables, and then using 
 $v_i(y) \leq cU_0(y)$ for $|y| \leq R_i$, we get
\begin{equation}\label{simples4}
\int_{|x| \leq s_i}u_i^{p_i} (x) dx = 
M_i^{-\frac{(p_i-1)n}{2}}M_i^{p_i}\int_{|y| \leq R_i}v_i^{p_i}(y)dy \leq 
CM_i^{-1}.
\end{equation}

On the other hand, by Claim 1,
\begin{eqnarray}\label{simples5}
\int_{s_i \leq |x| \leq \sigma}u_i^{p_i} (x) dx &\leq&
CM_i^{-\lambda_i p_i} 
\int_{s_i \leq |x| \leq \sigma}|x|^{(2-n+\delta)p_i}dx \nonumber\\
&\leq& CM_i^{-\lambda_i p_i}s_i^{(2-n+\delta)p_i+n}
\leq o(1)M_i^{-1}.
\end{eqnarray}

Claim 3 now follows from inequalities (\ref{simples4}) and (\ref{simples5}).

{\bf Claim 4:}
There exists $\sigma_1 > 0$ such that for all $0 < \sigma < \sigma_1$, there exists a constant
$C=C(\sigma)$ with, for every $i$,
$$
u_i(x_i)u_i(x) \leq C(\sigma)
$$
if $d(x,x_i) = \sigma$.

{\it Proof of Claim 4:}
From the discussion in Section \ref{confinvariance}, if we choose $\sigma_1 > 0$ small, we can suppose that $R_{g_i} \geq 0$.

Choose $0 < \sigma < \sigma_1$ small and define
$$
w_i(x) = u_i(x_{\sigma})^{-1}u_i(x),
$$
where $x_{\sigma}$ is  chosen so that $d(x_{\sigma},x_i)=\sigma$.
Note that
$$
L_{g_i}w_i + 
Ku_i(x_{\sigma})^{p_i-1} f_i^{-\tau_i}w_i^{p_i}=0.
$$

The Harnack inequality implies that, for every $\epsilon > 0$, 
there exists a constant $C_{\epsilon} > 0$ such that
$$
C_{\epsilon}^{-1} \leq w_i(x) \leq C_{\epsilon}
$$
if $d(x,\overline{x}) > \epsilon$. From Claim 1 we know that 
$u_i(x_{\sigma})^{p_i-1} \g 0$ as $i \g \infty$, and then standard elliptic
theory shows that, after maybe passing to a subsequence,
$$
w_i \g w \ {\rm in } \ 
C^2_{loc}(B_{\sigma}(\overline{x})),
$$
and $w$ satisfies
$$
L_{g}w = 0, w > 0.
$$

Since the blowup is isolated simple, the function $\hat{u}_i(r)$ is 
decreasing in the interval $(R_iM_i^{-\frac{p_i-1}{2}},\rho)$. Taking
the limit, we conclude that $\hat{w}(r)$ is decreasing in the whole
interval $(0, \rho)$. 

As a consequence, $w$ is singular  at
the origin. 

It follows from the results contained in the appendix in \cite{LIZHU99} that
\begin{equation}\label{simples6}
-\int_{B_{\eta}(x_i)}\Delta_{g_i}w_i = 
-\int_{\d B_{\eta}(x_i)} \frac{\d w_i}{\d \nu} = \\
-\int_{\d B_{\eta}(\overline{x})} \frac{\d w}{\d \nu}+o(1)> c > 0
\end{equation}
for each $i$, where $\eta > 0$ is sufficiently small.

On the other hand,
\begin{eqnarray}\label{simples7}
-\int_{B_{\eta}(x_i)}\Delta_{g_i}w_i &=&
\int_{B_{\eta}(x_i)}[Ku_i(x_{\sigma})^{-1} f_i^{-\tau_i}u_i^{p_i}-
c(n)R_{g_i}w_i] \nonumber\\
&\leq& K\int_{B_{\eta}(x_i)}u_i(x_{\sigma})^{-1} f_i^{-\tau_i}u_i^{p_i}
\leq cu_i(x_{\sigma})^{-1}M_i^{-1}.
\end{eqnarray}

Here we have used Claim 3. Claim 4 now follows from 
inequalities (\ref{simples6}) and (\ref{simples7}).

Now we are ready to prove Proposition \ref{propsimples}.

Suppose the inequality (\ref{estsimples}) is not true. Then there exists a sequence 
$\tilde{x}_i$, with $d(\tilde{x}_i,x_i) \leq \frac{\rho_1}{2}$, $\rho_1$ small, such that
\begin{equation}\label{simples8}
u_i(\tilde{x}_i)u_i(x_i)d(\tilde{x}_i,x_i)^{n-2} \g \infty
\end{equation}
as $i \g \infty$.

 Then, from Proposition \ref{propisol}, 
$R_iu_i(x_i)^{-\frac{p_i-1}{2}} \leq \tilde{r}_i \leq \frac{\rho_1}{2}$,
where $\tilde{r}_i=d(\tilde{x}_i,x_i)$.
Define
$$
\tilde{v}_i(y) 
= \tilde{r}_i^{\frac{2}{p_i-1}}u_i(\tilde{r}_iy), \ \ |y| < 2.
$$

Now it is not difficult to see that the origin is an isolated simple
blowup point for $\tilde{v}_i$, and Claim 4 implies, together with the Harnack inequality,
$$
\max_{|y|=1}\tilde{v}_i(0)\tilde{v}_i(y) \leq C.
$$
This contradicts the limit in (\ref{simples8}) and we finish the proof of inequality (\ref{estsimples}).

For the proof of inequality (\ref{estbaixo}), recall that the Green function always exists when $\rho_1$ is 
sufficiently small. 
Now observe that the inequality holds where 
$d_{g_i}(x,x_i) = R_iM_i^{-\frac{p_i-1}{2}}$ because $G_i(x_i,x)=O(r^{2-n})$ and where  $d_{g_i}(x,x_i) =   \rho_1$
because the Green function vanishes in this case and we are dealing with positive functions. Since 
$$
L_{g_i}(u_i(x_i)u_i) \leq 0 = L_{g_i}G_i,
$$
we can apply the maximum principle on the region $\{ x : R_iM_i^{-\frac{p_i-1}{2}} \leq d_{g_i}(x,x_i) \leq \rho_1 \}$
 to get the desired inequality.
\end{proof}


\begin{cor}\label{corsimples}
Under the hypotheses of Proposition \ref{propsimples}, after maybe passing to a subsequence, 
$$
M_i u_i(x) \g h \ {\rm in } \ 
C^2_{loc}(B_{\rho_1}(\overline{x}) \setminus \{\overline{x}\}),
$$
where $h$ is a positive solution to the linear equation $L_g(h)=0$,
with a nonremovable singularity at $\overline{x}$. (Here $g$ stands for
the limit metric.)
\end{cor}

\begin{proof}
Observe that the function $M_iu_i(x)$ satisfies
$$
L_{g_i}(M_iu_i(x)) + KM_i^{1-p_i}f_i^{-\tau_i}(x)(M_iu_i(x))^{p_i} = 0.
$$
The previous Proposition implies $M_iu_i(x)$ is uniformly bounded in 
compact sets contained in 
$B_{\rho_1}(\overline{x}) \setminus \{\overline{x}\}$, and then
standard elliptic estimates show that, after extracting a subsequence,
$M_iu_i(x) \g h$ in 
$C^2_{loc}(B_{\rho_1}(\overline{x}) \setminus \{\overline{x}\})$.
Since $M_i \g \infty$, $L_gh=0$.

Because of inequality (\ref{estbaixo}),  taking the limit,
one sees that  h
is singular. This finishes the proof.

\end{proof}


\section{Symmetry estimates}


In this section we will estimate the difference between solutions to our equation and standard symmetric functions,
which will be solutions to the corresponding critical equation in the Euclidean setting.

The next lemma gives us an estimate on $|v_i-U_0|$, depending on $M_i$ and $\tau_i=\frac{n+2}{n-2} - p_i$. This is the first step towards the symmetry
estimates, not depending on $\tau_i$, we will prove later.

\begin{lem}\label{lema1}
 Let $u_i$ be a sequence of positive functions satisfying equation (\ref{isol1}) and
$x_i \g \overline{x}$ be an isolated simple blowup point. Then there exists $\delta > 0$ such that
$$
|v_i(y)-U_0(y)| \leq C
\left\{
\begin{array}{lcr}
  \max \{M_i^{-2}, \tau_i\} &\ {\rm if } \ &n=4,5,\\
  \max \{(\log M_i)M_i^{-2}, \tau_i\} &\ {\rm if } \ &n=6,\\
  \max \{M_i^{-2}M_i^\frac{2(n-6)}{n-2}, \tau_i\} &\ {\rm if } \ &n\geq7,
\end{array}
\right.
$$
for $|y| \leq \delta M_i^{\frac{p_i-1}{2}}$, where $\tau_i = \frac{n+2}{n-2} - p_i$.
\end{lem}

\begin{proof}
Set $l_i = \delta M_i^{\frac{p_i-1}{2}}$, and 
 $\Lambda_i = \max_{|y|\leq l_i}|v_i - U_0| = v_i(y_i)-U_0(y_i)$, for a 
certain $|y_i| \leq l_i$.

We observe that if there exists a constant $c > 0$ such that $|y_i| \geq c l_i$ for every $i$, then inequality $v_i \leq cU_0$
 automatically 
implies the stronger
inequality 
$$
|(v_i-U_0)(y)| \leq C M_i^{-2},
$$
since 
$$
\Lambda_i = |v_i-U_0|(y_i) \leq C|y_i|^{2-n} \leq Cl_i^{2-n} \leq CM_i^{-2}.
$$
So for large $i$ we will have $|y_i| \leq \frac{l_i}{2}$. We are using that $M_i^{\tau_i} \g 1$ as $i \g \infty$.

Define 
$$ 
w_i(y) = \Lambda_i^{-1}(v_i(y)-U_0(y)).
$$ Then $w_i$ satisfies 
$$
L_{h_i}w_i + b_iw_i = Q_i(y),
$$
where 
$$
b_i(y) = K \tilde{f}_i^{-\tau_i} \frac{v_i^{p_i}-U_0^{p_i}}{v_i-U_0}(y),
$$
and
\begin{eqnarray}
&&Q_i(y) = \Lambda_i^{-1} \{
c(n)M_i^{-(p_i-1)}R_{g_i}(M_i^{-\frac{p_i-1}{2}}y)U_0(y) \\
&&+ M_i^{-(1+N)\frac{p_i-1}{2}}O(|y|^N)|y|(1+|y|^2)^{-\frac{n}{2}} +
K(U_0^\frac{n+2}{n-2}-\tilde{f}_i^{-\tau_i}U_0^{p_i}) \}\nonumber,
\end{eqnarray}
where $\tilde{f}_i(y) = f_i(M_i^{-\frac{p_i-1}{2}}y)$, $(h_i)_{kl}(y) = (g_i)_{kl}(M_i^{-\frac{p_i-1}{2}}y)$ and $O(|y|^N)$ comes from the expansion of the volume element in conformal normal coordinates and $N$ is as big as we want.

Since the blowup is isolated simple, from inequality $v_i \leq c U_0$,  it is easy to check, for example,
\begin{equation}\label{estbi}
b_i(y) \leq c(1+|y|)^{-3}
\end{equation}
for $|y| \leq l_i$.

We will choose $\delta$ small enough to guarantee the existence of the Green's function 
for the conformal Laplacian on a ball of radius $\delta$
, with respect to
a Dirichlet boundary condition.

The Green's representation formula gives 
\begin{equation}\label{green}
w_i(y) = \int_{B_i}G_{i,L}(y,\eta)(b_i(\eta)w_i(\eta) - Q_i(\eta))d\eta
- \int_{\d B_i}\frac{\d G_{i,L}}{\d \nu}(y,\eta)w_i(\eta)ds
\end{equation}
where $B_i$ stands for $B_{l_i}(0)$ and $G_{i,L}$ is the Green function of $L_{h_i}$ in 
$B_i$.

We will need the following lemma proved in  \cite{CHENLIN98}:


\begin{lem}\label{unicidade}
Suppose $w$ is a solution to the equation
\begin{equation}\label{eqw}
\Delta w + n(n+2)U_0^\frac{4}{n-2}w = 0 \ {\rm in} \ \Rn.
\end{equation}
 If $\lim_{|y| \g \infty} w(y) = 0$, then there exist constants $c_0,c_1,\dots,c_n$ such that
$$
w(y) = c_0(\frac{n-2}{2}U_0 + y \cdot \nabla U_0) + \sum_{j=1}^n c_j \frac{\d U_0}{\d y_j}.
$$
\end{lem}


{\bf Remark:}
The functions $\frac{n-2}{2}U_0 + y \cdot \nabla U_0$ and $\frac{\d U_0}{\d y_j}$, $j=1,\dots,n$, are solutions
to the equation (\ref{eqw}).

The proof of the Lemma \ref{lema1} is by contradiction.
Set
\begin{equation}\label{defti}
t_i = 
\left\{
\begin{array}{lcr}
  M_i^{-2}  &\ {\rm if } \ &n=4,5,\\
  (\log M_i)M_i^{-2} &\ {\rm if } \ &n=6,\\
  M_i^{-2}M_i^\frac{2(n-6)}{n-2} &\ {\rm if } \ &n\geq7.
\end{array}
\right.
\end{equation}

If the proposition is false, we necessarily have 
$$
\Lambda_i^{-1} \max \{t_i, \tau_i\} \g 0
$$
as $i \g \infty$, which  implies that
$$
\Lambda_i^{-1}t_i \g 0, \ \ \Lambda_i^{-1} \tau_i \g 0.
$$

Since $R = O(r^2)$ in conformal normal coordinates, we can get the following estimate:
\begin{eqnarray}\label{estqi}
|Q_i(y)| &\leq& c\Lambda_i^{-1} \{ M_i^{-\frac{8}{n-2}}|y|^2 (1+|y|)^{2-n}\nonumber\\
&+&M_i^{-(1+N)\frac{p_i-1}{2}}O(|y|^N)|y|(1+|y|^2)^{-\frac{n}{2}}\nonumber\\
&+&\tau_i (|\log U_0|+|\log \tilde{f}_i|) (1 + |y|)^{-n-2} \}.
\end{eqnarray}

Using the estimates (\ref{estbi}) and (\ref{estqi}),  we get from the Green's representation formula 
(\ref{green}) that $w_i$ is bounded in $C^2_{loc}$, and
\begin{equation}\label{estw45}
 |w_i(y)| \leq c[(1 + |y|)^{-1} + c\Lambda_i^{-1}t_i]
\end{equation}
for $|y| \leq \frac{\delta}{2}M_i^\frac{p_i-1}{2}$. We are using that $|w_i(y)| \leq C\Lambda_i^{-1}M_i^{-2}$
when $|y| = \frac{\delta}{2}M_i^\frac{p_i-1}{2}$, and also that  $|G_{i,L}(y,\eta)| \leq C|y-\eta|^{2-n}$
for $|y| \leq \frac{l_i}{2}$.

Then, by standard elliptic estimates, there 
exists a subsequence, also denoted $w_i$, converging to $w$ satisfying
$$
\left\{
\begin{array}{lr}
&\Delta w + n(n+2)U_0^{\frac{4}{n-2}}(y)w = 0 \ {\rm in }\  \Rn,\\
&|w(y)| \leq c(1+|y|)^{-1}.
\end{array}
\right.
$$

So, the Lemma \ref{unicidade} implies that
$$
w(y) = c_0(\frac{n-2}{2}U_0 + y \cdot \nabla U_0) + \sum_{j=1}^n c_j \frac{\d U_0}{\d y_j}.
$$
The conditions $w(0) = \frac{\d w}{\d y_j}(0) = 0$ show that $c_j=0$ for
every $j$, in other words, $w(y) \equiv 0$. From here we conclude that
$|y_i| \g \infty$ as $i \g \infty$.

This contradicts the estimate (\ref{estw45}) since $w_i(y_i)=1$ and $\Lambda_i^{-1}t_i \g 0$,
and this finishes the proof.
\end{proof}


In the next lemma, we estimate $\tau_i$. This result and the Lemma \ref{lema1} give us an estimate on
$|v_i-U_0|$ independent of $\tau_i$.

\begin{lem}\label{lema2}
Under the same hypotheses in Lemma \ref{lema1},
$$
\tau_i \leq C 
\left\{
\begin{array}{lcr}
M_i^{-2} &\ {\rm if }\ &n=4,5,\\
(\log M_i)M_i^{-2} &\ {\rm if }\  &n=6,\\
  M_i^{-2}M_i^\frac{2(n-6)}{n-2} &\ {\rm if }\  &n\geq7.
\end{array}
\right.
$$
\end{lem}

\begin{proof}
The proof will be again by contradiction and recall the definition (\ref{defti}). 
If the lemma is not true, then Lemma \ref{lema1} implies that  
$$
|v_i(y)-U_0(y)| \leq C \tau_i.
$$

Define
$$
w_i(y) = \tau_i^{-1}(v_i-U_0)(y),
$$
so $w_i$ is uniformly bounded. 
The equation satisfied by $w_i$ is
$$
L_{h_i}w_i + b_iw_i = \tilde{Q}_i(y),
$$
where 
$$
b_i(y) = K \tilde{f}_i^{-\tau_i}\frac{v_i^{p_i}-U_0^{p_i}}{v_i-U_0}(y)
$$
and
\begin{eqnarray}
\tilde{Q}_i(y) &=& \tau_i^{-1} \{
c(n)M_i^{-(p_i-1)}R_{g_i}(M_i^{-\frac{p-1}{2}}y)U_0(y) \nonumber\\
&+& M_i^{-(1+N)\frac{p_i-1}{2}}O(|y|^N)|y|(1+|y|^2)^{-\frac{n}{2}} \nonumber\\
&+& K(U_0^\frac{n+2}{n-2}-\tilde{f}_i^{-\tau_i}U_0^{p_i}) \}.
\end{eqnarray}

If the lemma  is not true, then $\tau_i^{-1}t_i \g 0$ as $i \g \infty$.

 We have
\begin{eqnarray}
|\tilde{Q}_i(y)| &\leq& c\tau_i^{-1} \{ M_i^{-\frac{8}{n-2}}|y|^2 (1+|y|)^{2-n} \nonumber\\
&+& M_i^{-(1+N)\frac{p_i-1}{2}}O(|y|^N)|y|(1+|y|^2)^{-\frac{n}{2}}\nonumber\\
&+&\tau_i (|\log U_0|+|\log \tilde{f}_i|) (1 + |y|)^{-n-2} \}.
\end{eqnarray}

By elliptic linear theory we can suppose $w_i \g w$ in compact subsets.

If $\psi(y) = \frac{n-2}{2}U_0(y) + y \cdot \nabla U_0(y)$, then, 
\begin{eqnarray}
&&\int_{|y| \leq \frac{l_i}{2}} \psi(y)\tau_i^{-1}(M_i^{-\frac{8}{n-2}}|y|^2 (1+|y|)^{2-n} \nonumber\\
&&\ \ \ \ + M_i^{-(1+N)\frac{p_i-1}{2}}O(|y|^N)|y|(1+|y|^2)^{-\frac{n}{2}})\g 0.
\end{eqnarray}
Note that when $i \g \infty$ we have:
$$
\tau_i^{-1}K(U_0^\frac{n+2}{n-2}-\tilde{f}_i^{-\tau_i}U_0^{p_i}) \g K (\log U_0(y)+\log f(\overline{x}))U_0^\frac{n+2}{n-2} 
$$
pointwise. It is not difficult to check that
$$
\int_{\Rn} \psi(y)U_0^\frac{n+2}{n-2}(y)dy = 0.
$$
Therefore we can conclude
$$
\lim_{i \g \infty} \int_{|y| \leq \frac{l_i}{2}} \psi(y)\tilde{Q}_i(y) dy =
n(n-2)\int_{\Rn} \psi(y)(\log U_0(y))U_0^\frac{n+2}{n-2}(y)dy.
$$

On the other hand, integration by parts shows that
\begin{eqnarray}
&&\int_{|y| \leq \frac{l_i}{2}} \psi(y)\tilde{Q}_i(y) dy = 
\int_{|y| \leq \frac{l_i}{2}} \psi(y)(L_{h_i} w_i + b_iw_i) dy  \nonumber\\
&&=\int_{|y| \leq \frac{l_i}{2}} (L_{h_i}\psi(y)+b_i\psi)w_i dy +
\int_{|y| = \frac{l_i}{2}} (\psi \frac{\d w_i}{\d r} - w_i \frac{\d \psi}{\d r}) d\sigma.
\end{eqnarray}

The integral on the boundary  goes to zero when $i \g \infty$ because
$$
\left\{
\begin{array}{lr}
& |\psi| = O(r^{2-n}),\ |\nabla \psi| = O(r^{1-n})  \\
& |w_i(\frac{l_i}{2})| \leq c\tau_i^{-1}M_i^{-2}, \ |\nabla w_i(\frac{l_i}{2})| \leq c\tau_i^{-1}M_i^{-2}l_i^{-1}. 
\end{array}
\right.
$$

Taking the limit when $i \g \infty$, we would have
$$
\lim_{i \g \infty} \int_{|y| \leq \frac{l_i}{2}} \psi(y)\tilde{Q}_i(y) dy =
\int_{\Rn} (\Delta \psi(y)+n(n+2)U_0^\frac{4}{n-2}\psi)w dy = 0
$$
because $\Delta \psi(y)+n(n+2)U_0^\frac{4}{n-2}\psi = 0$.

This is a contradiction because
\begin{equation}\label{eqlog}
n(n-2)\int_{\Rn} \psi(y)(\log U_0(y))U_0^\frac{n+2}{n-2}(y)dy > 0.
\end{equation}
To see this, first note that
$$
\psi(y) = \frac{n-2}{2}\frac{1-r^2}{(1+r^2)^\frac{n}{2}}.
$$
Observe that
\begin{eqnarray}\label{integral}
&&\int_{\Rn} \psi(y)(\log U_0(y))U_0^\frac{n+2}{n-2}(y)dy = \nonumber\\
&&-\frac{(n-2)^2}{4}\sigma_{n-1}\int_0^\infty \frac{1-r^2}{(1+r^2)^{n+1}}r^{n-1}\log (1+r^2)dr,
\end{eqnarray}
and after changing  variables $r=s^{-1}$, we get
$$
\int_0^{\infty}\frac{1-r^2}{(1+r^2)^{n+1}}r^{n-1}\log (1+r^2)dr = 2\int_1^\infty \frac{1-r^2}{(1+r^2)^{n+1}}r^{n-1}\log r dr.
$$
Now inequality (\ref{eqlog}) follows immediately, and that finishes the proof of the lemma.
\end{proof}


The Lemmas \ref{lema1} and \ref{lema2} together imply our symmetry estimate:

\begin{prop}\label{simetria}
Let $u_i$ be a sequence of positive functions satisfying equation (\ref{isol1}) and
$x_i \g \overline{x}$ be an isolated simple blowup point. Then there exists $\delta > 0$ such that
\begin{equation}\label{1symest}
|v_i(y)-U_0(y)| \leq C
\left\{
\begin{array}{lcr}
  M_i^{-2} &\ {\rm if } &\ n=4 \ {\rm or }\  5,\\
  (\log M_i)M_i^{-2} &\ {\rm if } &\ n=6,\\
  M_i^{-2}M_i^\frac{2(n-6)}{n-2} &\ {\rm if } &\ n\geq7
\end{array}
\right.
\end{equation}
for $|y| \leq \delta M_i^{\frac{p_i-1}{2}}$.
\end{prop}


When $n \geq 6$, by applying the same technique we can also get:

\begin{prop}\label{segsymest}
Let $u_i$ be a sequence of positive functions satisfying equation (\ref{isol1}) and
$x_i \g \overline{x}$ be an isolated simple blowup point.
Then
\begin{equation}\label{2symest}
|v_i(y)-U_0(y)| \leq C
\left\{
\begin{array}{lr}
 M_i^{-2}M_i^\frac{2}{n-2} (1+|y|)^{-1} &\ {\rm if } \ n=6 \\
 M_i^{-2}M_i^\frac{2(n-6)}{n-2} (1+|y|)^{6-n} &\ {\rm if } \ n\geq 7.
\end{array}
\right.
\end{equation}
\end{prop}

\begin{proof}

Set 
\begin{equation}
A_i = 
\left\{
\begin{array}{lr}
 M_i^{-2}M_i^\frac{2}{n-2}  &\ {\rm if } \ n=6 \\
 M_i^{-2}M_i^\frac{2(n-6)}{n-2}  &\ {\rm if } \ n\geq 7,
\end{array}
\right.
\end{equation}
and define
$$
w_i(y) = A_i^{-1}(v_i-U_0)(y)
$$
for $|y| \leq \delta M_i^{\frac{p_i-1}{2}}$. Then our previous proposition implies $w_i$ is uniformly bounded.
The equation satisfied is
$$
L_{h_i}w_i + b_iw_i = \tilde{Q}_i(y),
$$
where 
$$
b_i(y) = K \tilde{f}_i^{-\tau_i}\frac{v_i^{p_i}-U_0^{p_i}}{v_i-U_0}(y)
$$
and
\begin{eqnarray}
&&\tilde{Q}_i(y) = A_i^{-1} \{
c(n)M_i^{-(p-1)}R_g(M_i^{-\frac{p-1}{2}}y)U_0(y) \nonumber\\
 &&+ M_i^{-(1+N)\frac{p_i-1}{2}}O(|y|^N)|y|(1+|y|^2)^{-\frac{n}{2}}
+K(U_0^\frac{n+2}{n-2}-\tilde{f}_i^{-\tau_i}U_0^{p_i}) \}.
\end{eqnarray}
Then
\begin{eqnarray}
|\tilde{Q}_i(y)| &\leq& cA_i^{-1} \{ M_i^{-\frac{8}{n-2}}|y|^2 (1+|y|)^{2-n}\nonumber\\ 
&+& M_i^{-(1+N)\frac{p_i-1}{2}}O(|y|^N)|y|(1+|y|^2)^{-\frac{n}{2}}\nonumber\\ 
&+&\tau_i (|\log U_0|+|\log \tilde{f}_i|) (1 + |y|)^{-n-2} \}.
\end{eqnarray}
The Green's representation formula says that
\begin{equation}\label{green1}
w_i(y) = \int_{B_i}G_{i,L}(y,\eta)(b_i(\eta)w_i(\eta) - \tilde{Q}_i(\eta))d\eta
- \int_{\d B_i}\frac{\d G_{i,L}}{\d \nu}(y,\eta)w_i(\eta)ds,
\end{equation}
where $B_i$ stands for $B_{l_i}(0)$ and $G_{i,L}$ is the Green function of $L_{h_i}$ in 
$B_i$.

Since $|G_{i,L}(y,\eta)| \leq C|y-\eta|^{2-n}$, for $|y| \leq \frac{l_i}{2}$, we get
\begin{equation}\label{estwi2}
|w_i(y)| \leq c\{(1+|y|)^{-1} + cM_i^{-2}A_i^{-1} \} \leq c(1+|y|)^{-1}.
\end{equation}

If $n=6 \ {\rm or} \ 7$ the result follows multiplying the inequality (\ref{estwi2}) by $A_i$.
If $n \geq 8$ we plug the estimate (\ref{estwi2}) in the representation formula (\ref{green1})
until we reach
$$
|w_i(y)| \leq c(1+|y|)^{6-n}.
$$
Multiplying by $A_i$ we get the result.
\end{proof}

{\bf Remark 1:}
Once we have estimates (\ref{1symest}), (\ref{2symest}) on $v_i-U_0$, we can also get:
\begin{equation}
|\nabla(v_i-U_0)(y)| \leq C
\left\{
\begin{array}{lcr}
  M_i^{-2}(1+|y|)^{-1} & \ {\rm if }\  &n=4 \ {\rm or } \ 5,\\
  M_i^{-2}M_i^\frac{2}{n-2} (1+|y|)^{-2} &\ {\rm if } \ &n=6 \\
  M_i^{-2}M_i^\frac{2(n-6)}{n-2} (1+|y|)^{5-n}&\ {\rm if } \ &n\geq 7
\end{array}
\right.
\end{equation}
and
\begin{equation}
|\nabla^2(v_i-U_0)(y)| \leq C
\left\{
\begin{array}{lcr}
  M_i^{-2}(1+|y|)^{-2} &\ {\rm if } \ &n=4 \ {\rm or } \ 5,\\
  M_i^{-2}M_i^\frac{2}{n-2} (1+|y|)^{-3} &\ {\rm if } \ &n=6 \\
  M_i^{-2}M_i^\frac{2(n-6)}{n-2} (1+|y|)^{4-n}&\ {\rm if } \ &n\geq 7.
\end{array}
\right.
\end{equation}

{\bf Remark 2:}
If $h$ is as in the Corollary \ref{corsimples}, the estimates (\ref{1symest}) and (\ref{2symest}) imply that:
\begin{equation}
|h(x)-|x|^{2-n}| \leq C
\left\{
\begin{array}{lcr}
  1 &\ {\rm if } \ & n=4,5, \\
  |x|^{-1} &\ {\rm if } \ &n=6,\\
  |x|^{6-n} &\ {\rm if } \ &n \geq 7. 
\end{array}
\right.
\end{equation}
Since that gives the asymptotic behavior of the Green function of the conformal Laplacian in conformal normal
coordinates (see \cite{LEEPARKER87}), in some sense our symmetry estimates cannot be improved. 


\section{Local blowup analysis}


Now let us turn our attention to the applications of these symmetry estimates.

In the first application, we will show that the Weyl tensor of the metric has to vanish at an isolated  simple blowup point, when $n \geq 6$. This result had been proposed by R. Schoen (\cite{SCHOEN90}). This will allow us to use the Positive Mass Theorem in the proof of compactness of solutions to the Yamabe problem, when the dimension is 6 or 7.

\begin{thm}
Let $u_i$ be a sequence of positive functions satisfying equation (\ref{isol1}) and
$x_i \g \overline{x}$ be an isolated simple blowup point. If $n \geq 6$, then
$$
W_g(\overline{x}) = 0.
$$
\end{thm}

\begin{proof}

We will use the Pohozaev identity (\ref{pohozaev2}) to $u_i$ in  a ball of radius $r$:
\begin{eqnarray}\label{pohozaevweyl}
&&P(r,u_i)= \\
&&-\int_{|x| \leq r}(x^m\d_mu_i + \frac{n-2}{2}u_i)
((g^{kl}-\delta^{kl})\d_{kl}u_i + \d_kg^{kl}\d_lu_i)dx \nonumber\\
&&\ \ -c(n)\int_{|x| \leq r}(\frac12 x^k\d_kR + R)u_i^2 dx
+c(n)\frac{r}{2}\int_{|x|=r}Ru_i^2 d\sigma_r\nonumber\\ 
&&\ \ \ \ \ +(\frac{n}{p_i+1}-\frac{n-2}{2})\int_{|x| \leq r}Kf_i^{-\tau_i}u_i^{p_i+1}dx \nonumber\\
&&\ \ \ \ \ \ \ \ \ -\frac{\tau_i}{p_i+1}\int_{|x|\leq r}Kf_i^{-\tau_i-1}(x^m\d_mf_i)u_i^{p_i+1}dx.\nonumber
\end{eqnarray}
Using conformal normal coordinates we can get rid of the terms involving $|g|$.

Recall
\begin{eqnarray}
&&\ \ P(r,u_i) = \\
&&\int_{|x|=r}(\frac{n-2}{2}u_i\frac{\d u_i}{\d r}
-\frac{r}{2}|\nabla u_i|^2 + r|\frac{\d u_i}{\d r}|^2 + 
\frac{1}{p_i+1}Kf_i^{-\tau_i}ru_i^{p_i+1})d\sigma_r.\nonumber
\end{eqnarray}

Since we have that $M_i u_i \g h$ in the $C^2$ topology on compact subsets of $B_{\rho}(\overline{x}) \setminus \{\overline{x}\})$,
we conclude
\begin{equation}\label{estP}
M_i^2 |P(r,u_i)| \leq c < \infty.
\end{equation}

The same holds for
\begin{equation}\label{estR}
M_i^2|\int_{|x|=r}Ru_i^2 d\sigma(r)| \leq c < \infty.
\end{equation}
 The change of variables $y=M_i^\frac{p_i-1}{2}x$ yields
\begin{eqnarray*}
&&A_i(r) := \\
&&M_i^2 \{-\int_{|x| \leq r}(x^m\d_mu_i + \frac{n-2}{2}u_i)
((g^{kl}-\delta^{kl})\d_{kl}u_i + \d_kg^{kl}\d_lu_i)dx \nonumber\\
&&-c(n)\int_{|x| \leq r}(\frac12 x^k\d_kR + R)u_i^2 dx \} \nonumber\\
&&= -M_i^2 M_i^{2+(2-n)\frac{p_i-1}{2}}\int_{|y| \leq rM_i^\frac{p_i-1}{2}}
\{(y^m\d_mv_i 
+ \frac{n-2}{2}v_i)\\
&&((g^{kl}-\delta^{kl})(M_i^{-\frac{p_i-1}{2}}y)\d_{kl}v_i 
+ M_i^{-\frac{p_i-1}{2}}\d_kg^{kl}(M_i^{-\frac{p_i-1}{2}}y)\d_lv_i) \nonumber\\
&&+c(n)M_i^{-(p_i-1)}(\frac12 y^k\d_kR(M_i^{-\frac{p_i-1}{2}}y) 
+ R(M_i^{-\frac{p_i-1}{2}}y))v_i^2  \}dy.
\end{eqnarray*}

Note that $M_i^{2+(2-n)\frac{p_i-1}{2}} \g 1$.
Define
\begin{eqnarray}
&&\ \ \hat{A}_i(r) = \\
&& -M_i^2 M_i^{2+(2-n)\frac{p_i-1}{2}}  \int_{|y| \leq rM_i^\frac{p_i-1}{2}} \{(y^m\d_mU_0 
+ \frac{n-2}{2}U_0) \nonumber \\
&&((g^{kl}-\delta^{kl})(M_i^{-\frac{p_i-1}{2}}y)\d_{kl}U_0 
+ M_i^{-\frac{p_i-1}{2}}\d_kg^{kl}(M_i^{-\frac{p_i-1}{2}}y)\d_lU_0) \nonumber\\
&&+c(n)M_i^{-(p_i-1)}(\frac12 y^k\d_kR(M_i^{-\frac{p_i-1}{2}}y) 
+ R(M_i^{-\frac{p_i-1}{2}}y))U_0^2 \}dy.\nonumber
\end{eqnarray}

Then one can check
\begin{eqnarray}
&&\ \ |A_i(r) - \hat{A}_i(r)|  \\
&&\ \ \ \ \ \ \ \ \ \ \ \leq cM_i^2M_i^{-\frac{4}{n-2}}
\int_{|y| \leq rM_i^\frac{p_i-1}{2}} \{|v_i-U_0|(y)(1 + |y|)^{2-n} \nonumber \\
&&+|\nabla(v_i-U_0)|(1+|y|)^{3-n}
+ |\nabla^2(v_i-U_0)|(1+|y|)^{4-n}\}dy.\nonumber
\end{eqnarray}

If $n=6$, then
$$
|A_i(r) - \hat{A}_i(r)| \leq cM_i^{-\frac{2}{n-2}}\int_{|y| \leq rM_i^\frac{p_i-1}{2}}(1 + |y|)^{1-n}dy.
$$

When $n \geq 7$,
$$
|A_i(r) - \hat{A}_i(r)| \leq cM_i^{-\frac{4}{n-2}}M_i^{(n-6)\frac{2}{n-2}}
\int_{|y| \leq rM_i^\frac{p_i-1}{2}} (1 + |y|)^{8-2n}dy.
$$
Then
\begin{equation}\label{erroaichapeu}
|A_i(r) - \hat{A}_i(r)| \leq 
\left\{
\begin{array}{lcr}
  C \ &{\rm if } \ &n=6,7, \\
  C(\log M_i) \ &{\rm if } \ &n=8, \\
  CM_i^{(n-8)\frac{2}{n-2}} \ &{\rm if } \ &n\geq 9.
\end{array}
\right.
\end{equation}

Also observe that if we choose $r$ sufficiently small,
\begin{eqnarray}\label{tauisinal}
&&\ \ (\frac{n}{p_i+1}-\frac{n-2}{2})\int_{|x| \leq r}Ku_i^{p_i+1}dx  \\
&&-\frac{\tau_i}{p_i+1}\int_{|x|\leq r}Kf_i^{-\tau_i-1}(x^m\d_mf_i)u_i^{p_i+1}dx \geq 0,\nonumber
\end{eqnarray}
so we obtain, from identity (\ref{pohozaevweyl}), estimates (\ref{estP}), (\ref{estR}), (\ref{erroaichapeu})
and inequality (\ref{tauisinal}) that
$$
\hat{A}_i(r) \leq 
\left\{
\begin{array}{lcr}
  C \ &{\rm if } \ &n=6,7, \\
  C(\log M_i) \ &{\rm if } \ &n=8, \\
  CM_i^{(n-8)\frac{2}{n-2}} \ &{\rm if } \ &n\geq 9.
\end{array}
\right.
$$

We will need the Taylor series, for each $i$,
$$
R(x) = p_2(x) + p_3(x) + e(x)
$$
$$
(R + \frac12 x^k\d_kR)(x) = 2p_2(x) + \frac52 p_3(x) + e'(x)
$$
where $p_i$ is a homogeneous polynomial of degree $i$ and
$|e(x)|,|e'(x)| \leq c |x|^4$.

Let us denote by
$$
\tilde{u}_i(x) = M_i^{\frac{n-2}{4}\tau_i}\left( \frac{M_i^{-\frac{p_i-1}{2}}}{M_i^{-(p_i-1)} + |x|^2}\right)^\frac{n-2}{2}
$$
and it is not difficult to see that, changing variables, estimates on $|v_i-U_0|$ yield estimates on $|u_i-\tilde{u}_i|$.

Note that
$$
M_i^2\int_{|x| \leq r}(\frac12 x^k\d_ke + e)\tilde{u}_i^2 dx \leq
\left\{
\begin{array}{lcr}
C\ &{\rm if } \ &n=6,7,\\
C(\log M_i) \ &{\rm if } \ &n=8, \\
CM_i^{(n-8)\frac{2}{n-2}} \ &{\rm if } \ &n\geq 9.
\end{array}
\right.
$$

Then, using symmetry of $\tilde{u}_i$ 
$$
-c(n)\int_{|x| \leq r}2p_2(x)(M_i\tilde{u}_i)^2 dx \leq c
\left\{
\begin{array}{lcr}
  C \ &{\rm if } \ &n=6,7, \\
  C(\log M_i) \ &{\rm if } \ &n=8, \\
  CM_i^{(n-8)\frac{2}{n-2}} \ &{\rm if } \ &n\geq 9.
\end{array}
\right.
$$

But, on the other hand,  since  $\Delta R(0) = -\frac16 |W(0)|^2$ in conformal normal coordinates, we obtain
\begin{eqnarray*}
&&-c(n)M_i^2\int_{|x| \leq r}2p_2(x)\tilde{u}_i^2 dx \geq \\
&&\ \ \ \ \left\{
\begin{array}{lcr}
C|W(x_i)|^2M_i^2M_i^{-4\frac{2}{n-2}}(\log M_i) \ &{\rm if } \ &n=6,\\
C|W(x_i)|^2M_i^2M_i^{-4\frac{2}{n-2}} \ &{\rm if } \ &n \geq 7.
\end{array}
\right.
\end{eqnarray*}
So 
$$
|W_{g_i}(x_i)|^2 \leq 
\left\{
\begin{array}{lcr}
c(\log M_i)^{-1} &\ {\rm if } \ &n=6, \\
cM_i^{-\frac{2}{n-2}} &\ {\rm if } \ &n = 7, \\
cM_i^{-\frac{4}{n-2}}(\log M_i) &\ {\rm if } \ &n = 8, \\
cM_i^{-\frac{4}{n-2}} &\ {\rm if } \ &n \geq 9. 
\end{array}
\right.
$$
And taking the limit we are done.
\end{proof}


The next result concerns the local asymptotic analysis at a blowup point. It will be used together with the Positive Mass Theorem  to exclude the possibility of blowup phenomenon on manifolds not conformally diffeomorphic to the sphere.

Define
\begin{equation}\label{plinha}
P'(r,v) = 
\int_{|x|=r}(\frac{n-2}{2}v\frac{\d v}{\d \nu}
-\frac{r}{2}|\nabla v|^2 + r|\frac{\d v}{\d \nu}|^2)d\sigma(r).
\end{equation}

\begin{thm}\label{constneg}
Let $u_i$ be a sequence of positive functions satisfying equation (\ref{isol1}) and
$x_i \g \overline{x}$ be an isolated simple blowup point. and suppose $4\leq n \leq 7$. If $u_i(x_i)u_i \g h$ away from the origin, then
$$
\liminf_{r \g 0} P'(r,h) \geq 0.
$$
\end{thm}

\begin{proof}
That is another application of the Pohozaev identity (\ref{pohozaev})and the symmetry estimates:
\begin{eqnarray}
&& P(r,u_i)= 
-\int_{|x| \leq r}(x^m\d_mu_i + \frac{n-2}{2}u_i) \\
&&\ \ \ \ \ \ ((g^{kl}-\delta^{kl})\d_{kl}u_i + \d_kg^{kl}\d_lu_i 
-c(n)Ru_i)dx \nonumber \\
&&\ \ \ +(\frac{n}{p_i+1}-\frac{n-2}{2})\int_{|x| \leq r}Kf_i^{-\tau_i}u_i^{p_i+1}dx \nonumber\\
&&-\frac{\tau_i}{p_i+1}\int_{|x|\leq r}Kf_i^{-\tau_i-1}(x^m\d_mf_i)u_i^{p_i+1}dx.\nonumber
\end{eqnarray}

Firstly observe that
$$
M_i^2P(r,u_i) \g P'(r,h)
$$
as $i \g \infty$.

Secondly, as in the previous result, if $r$ is sufficiently small,
\begin{eqnarray}
&&(\frac{n}{p_i+1}-\frac{n-2}{2})\int_{|x| \leq r}Kf_i^{-\tau_i}u_i^{p_i+1}dx \nonumber\\
&&-\frac{\tau_i}{p_i+1}\int_{|x|\leq r}Kf_i^{-\tau_i-1}(x^m\d_mf_i)u_i^{p_i+1}dx \geq 0.
\end{eqnarray}

If 
\begin{eqnarray}
&&A_i(r) := M_i^2 \{-\int_{|x| \leq r}(x^m\d_mu_i + \frac{n-2}{2}u_i)
((g^{kl}-\delta^{kl})\d_{kl}u_i \nonumber\\
&&\ \ \ \ \ + \d_kg^{kl}\d_lu_i - c(n)Ru_i)dx 
 \}
\end{eqnarray}
and 
\begin{eqnarray}
&&\hat{A}_i(r) := M_i^2 \{-\int_{|x| \leq r}(x^m\d_m\tilde{u}_i + \frac{n-2}{2}\tilde{u}_i)
((g^{kl}-\delta^{kl})\d_{kl}\tilde{u}_i \nonumber\\
&&\ \ \ \ \ + \d_kg^{kl}\d_l\tilde{u}_i - c(n)R\tilde{u}_i)dx
 \},
\end{eqnarray}
then one can check
$$
|A_i(r)-\hat{A}_i(r)| \leq 
\left\{
\begin{array}{lcr}
cr^2 &\ {\rm if } \ &n=4,5, \\
cr   &\ {\rm if } \ &n=6,7.
\end{array}
\right.
$$

So
\begin{eqnarray}\label{limitep}
&&\liminf_{r \g 0} P'(r,h) \geq \\
&&c(n)M_i^2\liminf_{r \g 0}\int_{|x| \leq r}
(x^m\d_m\tilde{u}_i + \frac{n-2}{2}\tilde{u}_i)R\tilde{u}_idx .\nonumber
\end{eqnarray}

If $n=4$ or $5$, since $R = O(r^2)$, this automatically gives
$$
\liminf_{r \g 0} P'(r,h) \geq 0.
$$

If $n=6$ or $7$, the fact that $R = O(r^2)$ tells us that we only need to consider the second order term in the 
expansion of $R$. We are using the symmetry of $\tilde{u}_i$ to get rid of the third order term.

 The change of variables $y = M_i^\frac{p_i-1}{2}x$ implies that
\begin{eqnarray}
&&\ \ \ M_i^2\int_{|x| \leq r}
(x^m\d_m\tilde{u}_i + \frac{n-2}{2}\tilde{u}_i)R_{,ab}x^ax^b\tilde{u}_idx\\
&&=(1+o(1))M_i^\frac{2(n-6)}{n-2}
\int_{|y| \leq rM_i^\frac{p_i-1}{2}}(y^m\d_mU_0 + \frac{n-2}{2}U_0)R_{,ab}y^ay^bU_0dy.\nonumber
\end{eqnarray}

But
\begin{eqnarray}
&&\int_{|y| \leq rM_i^\frac{p_i-1}{2}}(y^m\d_mU_0 + \frac{n-2}{2}U_0)R_{,ab}y^ay^bU_0dy \nonumber\\
&&=\int_{0}^{rM_i^\frac{p_i-1}{2}}(r\d_rU_0 + \frac{n-2}{2}U_0)U_0r^{n+1} \frac{\Delta R}{n}\sigma_{n-1}dr\nonumber\\
&&=-\frac{1}{6n}\sigma_{n-1}|W(0)|^2\int_{0}^{rM_i^\frac{p_i-1}{2}}(r\d_rU_0 + \frac{n-2}{2}U_0)U_0r^{n+1}dr\nonumber\\
&&=-\frac{1}{6n}\sigma_{n-1}|W(0)|^2
(-\frac{n-2}{2}\int_{0}^{rM_i^\frac{p_i-1}{2}}\frac{r^{n+3}}{(1 + r^2)^{n-1}}dr \nonumber\\
&&\ \ \ \ \ \ \ \ \ \ \ \ \ \ \ \ \ \ \ \ \ \ \ \ \ \ \ 
+ \frac{n-2}{2} \int_{0}^{rM_i^\frac{p_i-1}{2}}\frac{r^{n+1}}{(1 + r^2)^{n-1}}dr). 
\end{eqnarray}

When $n=6$ and $i$ is large, this last expression is nonnegative because the first integral diverges while
the second one is finite. If $n=7$ and $i$ is large, again the expression is nonnegative since
\begin{equation}\label{relint}
\int_{0}^{\infty}\frac{r^{n+3}}{(1 + r^2)^{n-1}}dr = 
\frac{n+2}{n-6}\int_{0}^{\infty}\frac{r^{n+1}}{(1 + r^2)^{n-1}}dr.
\end{equation}

To see this, first observe that, if $m+1<2k$, integration by parts gives
$$
\int_{0}^{\infty}\frac{t^m}{(1+t^2)^k}dt = \frac{m-1}{2(k-1)}\int_{0}^{\infty}\frac{t^{m-2}}{(1+t^2)^{k-1}}dt.
$$
Since
$$
\int_{0}^{\infty}\frac{t^{m-2}}{(1+t^2)^{k-1}}dt = \int_{0}^{\infty}\frac{t^{m-2}}{(1+t^2)^k}dt + \int_{0}^{\infty}\frac{t^m}{(1+t^2)^k}dt,
$$
one gets
$$
\int_{0}^{\infty}\frac{t^m}{(1+t^2)^k}dt = \frac{m-1}{2k-m-1}\int_{0}^{\infty}\frac{t^{m-2}}{(1+t^2)^k}dt.
$$
Just choose $m=n+4$ and $k = n+1$ to obtain (\ref{relint}).

The proof is finished by using the inequality (\ref{limitep}).
\end{proof}


\section{Compactness theorem}

In this section we will prove the a priori estimates for the Yamabe problem in the non-locally conformally flat case,
for $4 \leq n \leq 7$.

The next Proposition is fundamental since it allows us to use the symmetry estimates we proved before.


\begin{prop}\label{isolatedsimple}
Suppose the blowup $x_i \g \overline{x}$ is isolated, $4 \leq n \leq 7$. Then it is also isolated simple.
\end{prop}

The proof of this statement  is just as in \cite{LIZHU99},
based on  Proposition \ref{constneg}.


Once we have established Proposition  \ref{isolatedsimple},
again following  \cite{LIZHU99}, we have:

\begin{thm}\label{blowupstructure2}
Suppose $u_i \in \mathbb{M}_{p_i}$ is a sequence satisfying 
$\max_M u_i \g \infty$, as $i \g \infty$. Then  $p_i \g \frac{n+2}{n-2}$.
Moreover, after passing to a subsequence,
\begin{enumerate}
\item the set $S= \{{\rm blowup \ points \ of } \ u_i \}$ is finite;
\item every blowup point of $u_i$ is an isolated simple blowup point.
\end{enumerate}
\end{thm}


Now let us turn to the statement and proof of the compactness theorem.

\begin{thm}
   Let $(M^n,g)$ be a smooth closed Riemannian manifold with positive Yamabe quotient, not conformally equivalent to $(\mathbb{S}^n, g_0)$. Assume $4 \leq n \leq 7$. Then, 
for every $\epsilon > 0$, there exists a 
positive constant $C = C(\epsilon, g)$ so that
$$
\left\{
\begin{array}{lr}
   &1/C \leq u \leq C \ {\rm and} \\
   &||u||_{C^{2,\alpha}(M)} \leq C 
\end{array}
\right.
$$
for every $u \in \cup_{1+\epsilon \leq p \leq \frac{n+2}{n-2}}\mathbb{M}_p$.
\end{thm}

\begin{proof}
Standard elliptic estimates and the Harnack inequality imply that it 
suffices to estimate $||u||_{C^0(M)}$.
Suppose, by contradiction, that $ \cup_{1+\epsilon \leq p \leq \frac{n+2}{n-2}}\mathbb{M}_p$ is not 
bounded in $C^0(M)$. This means that there exist $1+\epsilon \leq p_i \leq \frac{n+2}{n-2}$ and
$u_i \in \mathbb{M}_{p_i}$ with
$$
\max_M u_i \g \infty \ {\rm as } \ i \g \infty.
$$
From Theorem \ref{blowupstructure2} we know that this is only possible if $p_i \g \frac{n+2}{n-2}$.

Now Theorem \ref{blowupstructure2} implies that, after possibly passing to a subsequence, $u_i$ has N isolated simple
blow-up points $x_i^{(1)} \g x^{(1)}, \dots, x_i^{(N)} \g x^{(N)}$, for some 
integer N.

Define $w_i(x) = u_i(x_i^{(1)})u_i(x)$. We can suppose, for example, that
$$
u_i(x_i^{(1)}) = \min \{u_i(x_i^{(1)}),\dots,u_i(x_i^{(N)})\} \ \ {\rm for \ all} \ \ i.
$$
Proposition \ref{propsimples} then implies that there exists $\rho,c > 0$ such that
\begin{equation}\label{comp2}
w_i(x) \leq cd(x,x_i^{(j)})^{2-n} \ \ {\rm when } \ \ d(x,x_i^{(j)})\leq \rho, \ 
1\leq j \leq N.
\end{equation}

On the other hand, we know that the sequence $u_i$ is uniformly bounded in
$M \setminus \cup_{j=1}^{j=N}B_{\frac{\rho}{4}}(x^{(j)})$, since there is no blowup point in that region.  Then 
the Harnack inequality implies $w_i$ is uniformly bounded in 
$M \setminus \cup_{j=1}^{j=N}B_{\frac{\rho}{2}}(x^{(j)})$. This and 
inequality (\ref{comp2}) imply, after passing to a subsequence,
$$
u_i(x_i^{(1)})u_i(x) \g h(x) = \sum_{j=1}^{N}a_jG_{x^{(j)}}(x) + b(x) 
$$ 
in $C^2_{loc}(M \setminus \{x^{(1)}, \dots, x^{(N)}\})$,
where $a_1,\dots,a_N$ are nonnegative  constants, $G_{x^{(j)}}$ is the Green function of the 
conformal Laplacian with pole at $x^{(j)}$ and $b(x)$ is a regular $C^2$ function satisfying
$L_g (b) = 0$ in M.
Note that, since the first eigenvalue of the conformal Laplacian $L_g$ is positive, $b \equiv 0$.

Since the Green functions considered are positive and $a_1 > 0$ because of 
the Corollary \ref{corsimples}, we obtain the following expansion (\cite{LEEPARKER87}) :
$$
h(x) = ad(x,x^{(1)})^{2-n} + A + O(r), \ {\rm for } \ 4 \leq n \leq 6
$$
or
$$
h(x) = ad(x,x^{(1)})^{2-n} - cR_{;ij}x_ix_jr^{-3} + A + O(r), \ {\rm for } \ n=7.
$$
 The Positive Mass Theorem asserts that $A>0$.

We want to compute the limit 
$$
\liminf_{r \g 0}P'(r,h).
$$ 

Since $r^{2-n}$ is harmonic, definition (\ref{plinha}) gives
\begin{eqnarray}
&P'(r,h)= \int_{|x|=r} \{-\frac{(n-2)^2}{2}aAr^{1-n}+O(r^{2-n})\} d\sigma_r.
\end{eqnarray}
When $n=7$, we  use symmetry and the fact that $\Delta R(0)=-\frac16 |W(0)|^2 = 0$.

 Then 
$$
\liminf_{r \g 0}P'(r,h) = -\frac{(n-2)^2}{2}aA\sigma_{n-1} < 0,
$$
by the Positive Mass Theorem.

This contradicts Theorem \ref{constneg} and it finishes the proof.

\end{proof}

\footnotesize

{\sc Fernando C. Marques, Instituto de Matem\'{a}tica Pura e Aplicada (IMPA), Estrada Dona Castorina 110, 22460-320, Rio de Janeiro - RJ,
Brazil}

{\it Email address : \ } coda@impa.br


\bibliography{fernando.bib}

\begin{thebibliography}{10}

\bibitem{AUBIN76}
T.~Aubin.
\newblock \'{E}quations diff\'erentielles non lin\'eaires et probl\'eme de
  {Y}amabe concernant la courbure scalaire.
\newblock {\em J. Math. Pures Appl.}, 55:269--296, 1976.

\bibitem{CAFFAGIDASSPRUCK89}
L.~Caffarelli, B.~Gidas, and J.~Spruck.
\newblock Asymptotic symmetry and local behavior of semilinear elliptic
  equations with critical {S}obolev growth.
\newblock {\em Comm. Pure Appl. Math}, 42:271--297, 1989.

\bibitem{CHENLIN98}
C.~C. Chen and C.~S. Lin.
\newblock Estimates of the scalar curvature equation via the method of moving
  planes. ii.
\newblock {\em J. Diff. Geom.}, 49:115--178, 1998.

\bibitem{DRUET2003}
O.~Druet.
\newblock From one bubble to several bubbles: the low-dimensional case.
\newblock {\em J. Diff. Geom.}, 63(3):399--473, 2003.

\bibitem{DRUET2004}
O.~Druet.
\newblock Compactness for {Y}amabe metrics in low dimensions.
\newblock {\em Int. Math. Res. Not.}, 23:1143--1191, 2004.

\bibitem{HANLI99}
Z.~C. Han and Y.~Y. Li.
\newblock The {Y}amabe problem on manifolds with boundary: existence and
  compactness results.
\newblock {\em Duke Math. Journal}, 99(3):489--542, 1999.

\bibitem{LEEPARKER87}
J.~Lee and T.~Parker.
\newblock The {Y}amabe problem.
\newblock {\em Bull. Amer. Math. Soc.}, 17:37--91, 1987.

\bibitem{LI95}
Y.~Y. Li.
\newblock Prescribing scalar curvature on $\mathbb{S}^n$ and related problems,
  part {I}.
\newblock {\em J. Differential Equations}, 120:319--410, 1995.

\bibitem{LIZHU99}
Y.~Y. Li and M.~Zhu.
\newblock Yamabe type equations on three dimensional {R}iemannian manifolds.
\newblock {\em Comm. in Contemporary Mathematics}, 1(1):1--50, 1999.

\bibitem{MARQUES2003}
F.~C. Marques.
\newblock {\em Existence and compactness theorems on conformal deformation of
  metrics}.
\newblock PhD thesis, Cornell University, 2003.

\bibitem{OBATA71}
M.~Obata.
\newblock The conjectures on conformal transformations of {R}iemannian
  manifolds.
\newblock {\em J. Diff. Geom.}, 6:247--258, 1971.

\bibitem{POLLACK93}
D.~Pollack.
\newblock Compactness results for complete metrics of constant positive scalar
  curvature on subdomains of $\mathbb{S} \sp n$.
\newblock {\em Indiana Univ. Math. J.}, 42(4):1441--1456, 1993.

\bibitem{SCHOEN91}
R.~Schoen.
\newblock On the number of constant scalar curvature metrics in a conformal
  class.
\newblock In {\em Differential Geometry: A Symposium in Honor of Manfredo do
  Carmo}.

\bibitem{SCHOEN90}
R.~Schoen.
\newblock A report on some recent progress on nonlinear problems in geometry.
\newblock In {\em Surveys in differential geometry (Cambridge, MA, 1990)},
  pages 201--241. Lehigh Univ., Bethlehem, PA.

\bibitem{SCHOEN84}
R.~Schoen.
\newblock Conformal deformation of a {R}iemannian metric to constant scalar
  curvature.
\newblock {\em J. Differential Geometry}, 20:479--495, 1984.

\bibitem{SCHOENYAU79}
R.~Schoen and S.~T. Yau.
\newblock On the proof of the positive mass conjecture in {G}eneral
  {R}elativity.
\newblock {\em Comm. Math. Phys.}, 65:45--76, 1979.

\bibitem{SCHOENZHANG96}
R.~Schoen and D.~Zhang.
\newblock Prescribed scalar curvature on the n-sphere.
\newblock {\em Calc. Var. and PDEs}, 4:1--25, 1996.

\bibitem{TRUDINGER68}
N.~Trudinger.
\newblock Remarks concerning the conformal deformation of a {R}iemannian
  structure on compact manifolds.
\newblock {\em Ann. Scuola Norm. Sup. Pisa Cl. Sci.}, 22(3):165--274, 1968.

\bibitem{YAMABE60}
H.~Yamabe.
\newblock On a deformation of {R}iemannian structures on compact manifolds.
\newblock {\em Osaka Math. J.}, 12:21--37, 1960.

\end{thebibliography}

\bibliographystyle{plain}

\end{document}